\documentclass[final]{aupl}

\usepackage{
amsfonts,
amsmath,
latexsym,
amssymb,
amscd,
enumerate,
verbatim,
mathrsfs,
graphicx,
centernot,
accents,
color,}

\newcommand{\labbel}[1]{\label{#1} [[{\bf #1}]]}  
\renewcommand{\labbel}{\label}

\newtheorem{theorem}{Theorem}[section]
\newtheorem{lemma}[theorem]{Lemma}

\newtheorem{corollary}[theorem]{Corollary}

\newtheorem{void}[theorem]{} 
\newtheorem*{claim*}{Claim}

\newtheorem*{theorem*}{Theorem}
\newtheorem*{proposition*}{Proposition}
\newtheorem*{corollary*}{Corollary}
\newtheorem*{lemma*}{Lemma}
\newtheorem*{scholion*}{Scholion}

\theoremstyle{definition}

\theoremstyle{remark}
\newtheorem{remark}[theorem]{Remark}

\newtheorem*{remark*}{Remark}
\newtheorem*{remarks*}{Remarks}

\newtheorem*{observation*}{Observation}

\newcommand{\brfrt}{\hspace{0 pt}}
 
 \allowdisplaybreaks[1]

\numberwithin{equation}{section}

\begin{document}

\title{A note on non-generators in partially ordered sets}

\author{Paolo Lipparini}

\urladdr{http://www.mat.uniroma2.it/~lipparin}

\address{Dipartimento di Matematica\\Viale della  Ricerca
Non Generata\\Universit\`a di Roma ``Tor Vergata'' 
\\I-00133 ROME ITALY}



\date{\today}

\thanks{
Work performed under the auspices of G.N.S.A.G.A. Work 
partially supported by PRIN 2012 ``Logica, Modelli e Insiemi''.
The author acknowledges the MIUR Department Project awarded to the
Department of Mathematics, University of Rome Tor Vergata, CUP
E83C18000100006.}

\begin{abstract}
A folklore argument shows that 
Frattini's characterization of non-\brfrt generators
works in the framework of algebraic partially ordered sets.
We provide characterizations of non-generators
in arbitrary partially ordered sets.
The validity of some characterizations is equivalent
to  Zorn's Lemma, hence to the Axiom of choice.
We notice that working on closure spaces or posets 
provides no essential improvement. 
\end{abstract}

\maketitle  

\section{Introduction} \labbel{intro} 

Most results in the present note should be considered
as folklore, but are possibly new at the
level of generality we shall introduce.
In any case, we hope that the note might help
in clarifying some issues which sometimes appear
obscure or neglected.

In a classical paper  
G. Frattini \cite{Fr} showed that
an element $g$ of some  group $\mathbf G$ is a non-generator
if and only if  $g$ belongs to all the maximal subgroups of 
$\mathbf G$. 
Motivated by \cite{Fr}, the intersection of all
the maximal subgroups of a group is called the \emph{Frattini}
subgroup. 
Many years later N. Jacobson \cite{Ja} 
introduced a similar notion in ring theory.
Scattered generalizations to other kinds of 
structures appeared in the literature,
though not  as successful as Frattini's and Jacobson's notions, so far.
See, e.~g., \cite{ADS,BS,KST,MSS,Sc,S,T} 
and further references in the quoted sources. 
In a recent illuminating paper G. Janelidze
 \cite{Jz} 
exploited Frattini's methods
in a categorical context.

With a more modest aim, here we try 
to make the arguments work in the barest possible
setting, that is, partially ordered sets, \emph{posets}, for short,  with no
finiteness, completeness or algebraic assumption.
This should be explained a bit, since
Frattini's work deals with the notion of ``generation'',
while the meaning of ``generate'' is somewhat ambiguous 
within the realm of
posets.   
Indeed, at first sight, the simplest way to insert 
Frattini's characterization into an  abstract framework
would be to consider the power set 
$\mathcal P(G)$ of some set  $G$,
together with a notion of ``generation'',
``hull'', or ``closure''.
We shall see that 
Frattini's characterization works in a much more
general framework. Of course, in the specific case
of groups,   if $\mathbf G$ is a group
and $X \subseteq G$, the  ``closure''  $\langle X \rangle$
of $X$ 
is interpreted as the subgroup of  $\mathbf G$
generated by $X$. 

Abstracting from the above example,
we obtain the notion of a \emph{closure space},
roughly, a structure like a topological space, 
except that the union of two ``closed'' subsets
is not necessarily assumed to be closed. 
Closure spaces have 
a long history and have found an unexpected number of
applications in extremely disparate fields \cite{E}.
Just as an example, the above starting example can be generalized to 
every algebraic structure.
If $\mathbf A$ is an arbitrary algebraic structure, 
then $\mathcal P(A)$ becomes a closure space
when the (domains of) subalgebras of $\mathbf A$ 
 are considered as closed subsets,
 possibly including the empty set.
However, as we mentioned, there are many more possible examples,
we refer again to \cite{E} for more information.
  
It has been noticed in \cite{Jz} that, as far as Frattini's characterization is concerned,
we can deal just with the ``closed'' subsets,
the lattice of subgroups, in the original example.
Indeed, an element $g$ of some group is a non-generator   
if and only if the subset $\langle g \rangle$ is a non-generator in the 
lattice of subgroups of $\mathbf G$. 
Actually, we do not need the full lattice structure on the set
of all subgroups of a group, it is enough to
deal with the join-semilattice structure.
The same remark applies to an arbitrary closure space;
in fact, to any poset endowed with a closure operation.
Namely, rather than dealing with the full closure poset,
it is enough to deal with the lattice of closed elements.

Here we go just one step further. 
In order to exploit Frattini's arguments,
we do not need the actual notion or join or  ``generation'',
we only need to know what ``generating the whole structure''
means. Hence we can work in an arbitrary poset $\mathbf P$, not necessarily
a join-semilattice, provided that $\mathbf P$
has a maximum element $1$.  
A subset $X \subseteq P$ is meant to ``generate'' $1$
if there is no $p \in P$ such that    $p < 1$ and  $x \leq p$, for every 
$x \in X$. As we shall mention, the assumption of the existence 
of a maximum is just for convenience, since all the motivating 
examples have a maximum.

However---already in the setting of complete lattices---the
generalization of Frattini's characterization does not
necessarily hold, unless some algebraicity assumption is made
\cite[Remark 2.3]{Jz}. 
We show that algebraicity is enough to guarantee
that Frattini's characterization holds in the realm of
posets; moreover, we weaken it to a form of algebraicity
relative only to ``unbounded limits''. 
Actually, we notice that we only need 
algebraicity with respect to chains all whose members
are not extendable to a maximal element.
Furthermore, we give a characterization of non-generators in a poset
even when no algebraicity assumption is made.

In a final section we notice that, as far as Frattini's
characterization is concerned,
dealing with closure spaces or closure posets
provides no essential gain.
All the significant results seem to follow
already from the results provable for posets
with no further structure.

\section{Non-generators in partially ordered sets} \labbel{ngsec} 

Henceforth,  $\mathbf P=(P, \leq)$
is a fixed  \emph{partially ordered set},
 \emph{poset}, for short.
For convenience, we suppose that $\mathbf P$
has a maximum $1$.   
An element $a \in P$
is a \emph{non-generator}
if $a \vee e = 1$
implies $e=1$,
for every $e \in P$.
Notice that we do not assume that $P$
is a semilattice; the meaning of 
$a \vee e = 1$ is that there is no 
$p \in P$ such that $p < 1$ and  both 
$a \leq p$ and $e \leq p$.
In other words,  $a \vee e = 1$
means that $a$ and  $e$ have no common bound
 $<1$. 
The notation is consistent
with the standard notion of join
when  considered as a partial operation.

An element $p \in P$
is \emph{maximal$_{_{<1}}$}
if $p < 1$ and there is no $q \in P$
such that  $p < q < 1$.
Thus $p $
is maximal$_{_{<1}}$
in $\mathbf P$
if and only if 
$p$ is maximal in the standard sense
in the poset induced by $\mathbf P$
on $P \setminus \{ 1 \}  $.
Alternative expressions
for maximal$_{_{<1}}$
are \emph{proper maximal} or
 \emph{almost maximal}.
However, the use of maximal 
for what we call maximal$_{_{<1}}$
is quite standard, though ambiguous.
We have preferred to use a terminology
as closed as possible to the standard one,
just adding a subscript in order to avoid ambiguity.

Similarly to the finite case, if $X \subseteq P$,
then $\bigvee X=1$
means that there is no  $p \in P$ such that $p < 1$ and  
$x < p$, for every $x \in X$.  
Intuitively, $\bigvee X=1$ should always be intended
as ``$X$ generates the whole structure'',
hence $\langle X \rangle= 1$
could be an alternative notation.  
We have favored the notation compatible
with the standard notation used in (partial) lattice theory.

We could  rephrase all the results here
in the case when $\mathbf P$
has no maximum; in this case, 
 $a \vee b = 1$ should be replaced by
``$a$ and  $b$ have no common bound'', and similarly for 
$\langle X \rangle= 1$. 
 The assumption of the existence of $1$
is just because all the concrete examples have indeed  a maximum.

 \begin{lemma} \labbel{join}
Let $\mathbf P$ be a poset with maximum $1$
and let $a, b \in P$. 
  
(i)
If $a $ is a non-generator
and  $b \leq a$,
then $b$ is a non-generator.

(ii)
If $a , b $ are non-generators
and the join $a \vee b$ exists in $\mathbf P$,
then $a \vee b$ is a non-generator.

(iii)
If $a $ is a non-generator,
$X \subseteq P$ and 
$\bigvee (X \cup \{ a \} )=1$,
then   $\bigvee X =1$.
 \end{lemma}

 \begin{proof}
(i) If $ b\vee e=1$, for some $e \in P$, then
$ a\vee e=1$, hence $e=1$, since $a$ is a non-generator.

(ii)
Suppose that 
$(a \vee b) \vee e =1$;
this means that 
there is no 
$p \in P$ such that $p < 1$,
$a \leq p$, $b \leq p$ and $e \leq p$.

Suppose that $b \leq d$ and 
$e \leq d$, for some $d$.
Then $a \vee d= 1$, since otherwise
there is some $p < 1$ such that $a \leq p$
and $d \leq p$, thus also 
$b \leq p$ and $e \leq p$, a contradiction.
Since $a$ is a non-generator, then
$d=1$.

We have showed that there is no $d<1$
such that $b \leq d$ and 
$e \leq d$, that is $b \vee e =1$.     
Since $b$ is a non-generator, then
$e=1$, what we had to show.  

(iii)
Suppose by contradiction that $\bigvee X =1$
fails. Then there is $r \in P$ such that $r<1$
and $x \leq r$, for every $x \in X$.
Then $r \vee a=1$, since     $\bigvee (X \cup \{ a \} )=1$.
Since $a$ is a non-generator, then
$r=1$, a contradiction.   
\end{proof}    

For short, (i) and (ii) entail that in a join semilattice the set of
non-generators is an ideal.
However, Lemma \ref{join}
holds even for posets which are not semilattices. 
Clause (iii) means that being  a non-generator with respect to 
elements is the same as being  a  
non-generator with respect to subsets.

As we mentioned, G. Frattini showed that 
in a group the set of non-generators
is the intersection of all the maximal 
(maximal$_{_{<1}}$, in the present terminology) subgroups.
For arbitrary posets, it is not
always the case that some element
$a$ is a non-generator if and only if 
$a \leq p$, for every maximal$_{_{<1}}$  $p$.

A \emph{chain} in a poset 
is a subset which is 
linearly ordered by the induced order.
We shall always assume that 
chains are nonempty.
By a slight abuse of notation,
we shall denote chains as sequences.
In what follows, we could always assume that chains are well-ordered,
however we do not give details, since this remark shall never
be used.
 
If $\mathbf Q$ is a poset, 
we say that a chain $( e_i) _{i \in I} $ in $\mathbf Q$
is \emph{unbounded}
if there is no $q \in Q$  
such that $e_i \leq q$, for every $i \in I$.

We rephrase Zorn's Lemma as follows.

\begin{void} \labbel{zl}
\textsc{(Zorn's Lemma)} For every nonempty poset $\mathbf Q$,
either $\mathbf Q$ has a maximal element, 
or there is an unbounded chain in  
$\mathbf Q$.
 \end{void}

Suppose that $\mathbf P$ is a poset with maximum $1$. 
We say that a chain $( e_i) _{i \in I} $ in $\mathbf P$
is \emph{unbounded$_{_{<1}}$} 
if
$e_i <1$, for every $i \in I$, and moreover 
 there is no $p \in P$  
such that $p<1$ and $e_i \leq p$, for every $i \in I$.
Thus, for every poset $\mathbf P$
with maximum,  a chain 
 is unbounded in $P \setminus \{ 1 \} $ if and only if 
it 
is unbounded$_{_{<1}}$  in $\mathbf P$
and $1$ does not belong to the chain.  

\begin{theorem} \labbel{cpp}
Suppose that $\mathbf P$ is a poset with maximum $1$.
If $a \in P$, the following conditions are equivalent.
 \begin{enumerate}[(1)]  
 \item 
The element $a$ is a non-generator.
\item
Both (i) $a \leq p$, for every maximal$_{_{<1}}$ $p$, and
  \begin{enumerate}    
\item[(ii)]
for every $e \in P$ such that $e <1$ and
$e$  cannot be extended to a maximal$_{_{<1}}$  element,
there is $r \in P$ such that $r<1$, $a < r$ and $e < r$.      
 \end{enumerate}   
\item
Both (i) $a \leq p$, for every maximal$_{_{<1}}$ $p$, and
 \begin{enumerate}    
\item[(iii)]
for every unbounded$_{_{<1}}$
chain $( e_i) _{i \in I} $ and for every 
$i \in I$ (equivalently, for some $i \in I$),
there is $r \in P$ such that 
$a \leq r$ and $e_i \leq r$.   
 \end{enumerate}   
\item
Both (i) $a \leq p$, for every maximal$_{_{<1}}$ $p$, and
 \begin{enumerate}    
\item[(iv)]
whenever  $( e_i) _{i \in I} $
 is an unbounded$_{_{<1}}$
chain such that no element of the chain
can be extended to a maximal$_{_{<1}}$ element,
 then, for every 
$i \in I$ (equivalently, for some $i \in I$),
there is $r \in P$ such that 
$a \leq r$ and $e_i \leq r$.   
 \end{enumerate}   
\end{enumerate}  
\end{theorem} 

 \begin{proof}
(1) $\Rightarrow $  (2)(i) If $p$ is maximal
and not $a \leq p$,
then $a \vee p= 1$.
Since $a$ is a non-generator, then
$p=1$, contradicting $p<1$ from the definition of
maximality.

(1) $\Rightarrow $  (2)(ii) is trivial; actually,
if $a$ is a non-generator,  
then \emph{for every} $e \in P$, if  $e <1$, then
there is $r \in P$ such that $r<1$, $a < r$ and $e < r$.      

(2) $\Rightarrow $  (1) Assume (2) and suppose that $e \in P$ and  
$a \vee e=1$. We have to show that $e=1$.
Suppose not. If $e$ cannot be extended to a maximal$_{_{<1}}$  element,
then (ii) contradicts $a \vee e=1$.
Otherwise, there is some  maximal$_{_{<1}}$ $p$
such that $e \leq p$. But then   
$a \vee e=1$ implies $a \vee p=1$,
which  is impossible, since $a \leq p$, by (i). 
We have reached a contradiction in each case, hence 
$e=1$.

(1) $\Rightarrow $  (3)(iii)
is immediate since, by construction,
 $e_i < 1$, for every  $i \in I$.

(3) $\Rightarrow $  (4) is  trivial.

(4) $\Rightarrow $  (1) 
Assume (4) and let $e <1$.
We want to show that there is $r<1$
such that $a  \leq r$ and  $e \leq r$.  
Apply Zorn's Lemma \ref{zl}
to the poset $Q=\{ \,  b \in P  \mid e \leq b < 1  \,\}$,
with the order induced by $\mathbf P$.  
Notice that $Q$ is nonempty, since $e <1$. 
 If $Q$ has a maximal element $p$,
then $p$ is maximal$_{_{<1}}$ in $\mathbf P$,
hence $a \leq p$ by (i) and $e \leq p$ by construction, so we can take
$r=p$.
Otherwise, $\mathbf Q$ has an unbounded chain,
which henceforth is an unbounded$_{_{<1}}$ chain
in  $\mathbf P$. Then $r$ is given by (iv).  
\end{proof} 

While the equivalence of (1) and (2)
in \ref{cpp} is trivial, though useful,
the other equivalences are less trivial;
actually they imply back Zorn's Lemma.
The following corollary is stated in some set theory
in which the Axiom of choice is not assumed, say,
Zermelo Fraenkel set theory (ZF).

\begin{corollary} \labbel{z}
(ZF) The assertion that Clause
(3) in Theorem \ref{cpp} implies Clause (1)
is equivalent to Zorn's Lemma, hence also to the Axiom of choice.
 \end{corollary}

 \begin{proof} 
Suppose that Clause (3) implies Clause (1),
we shall prove Zorn's Lemma.

Let $\mathbf Q$ be a nonempty poset.
If $\mathbf Q$ has some maximal element, we are done.
Otherwise, let $P= Q \cup \{ 1, a \} $
and order $P$ by adding the conditions
$a <1$, 
$q < 1$, for every $ q \in Q$
and let $a$ be incomparable with every element of $Q$.
Since $Q$ is nonempty, then $a$  fails to be  a non-generator,
hence Clause (3) must fail for $a$.
Since $\mathbf Q$ has no maximal element, the only 
maximal element of $\mathbf P$ is $a$, hence (3)(i) is satisfied.
Thus (3)(ii) fails, $\mathbf P$ has some unbounded$_{_{<1}}$
chain, hence  $\mathbf Q$ has some unbounded
chain.

All the rest is trivial or well-known.
\end{proof}

An element $a$ of some poset $\mathbf P$  with a maximum
$1$ 
is \emph{$1$-compact} if, whenever
$( e_i) _{i \in I} $  is a chain of elements of $P$
and $\bigvee _{i \in I} e_{i}=1$,
then $a \leq e_i$, for some $i \in I$.  
As a standard argument, Theorem  \ref{cpp} 
implies that a $1$-compact element $a$ is a non-generator
if and only if  $a$ is contained in every
maximal$_{_{<1}}$ element. 

In fact, the argument works
for an even larger class of elements.
Say that $a$ 
is \emph{weakly-$1$-compact} if, whenever
$( e_i) _{i \in I} $  is a chain of elements of $P$
and $\bigvee _{i \in I} e_{i}=1$,
then there are 
some $i \in I$ and some $r \in P$ 
such that $r<1$, $ e_i \leq r$ and $a \leq r$. 
As an even weaker notion,
say that
$a$ 
is \emph{very-weakly-$1$-compact}
if the above condition applies, limited
to those chains 
$( e_i) _{i \in I} $ such that no 
$e_i$ can be extended to a maximal$_{_{<1}}$  element.

\begin{corollary} \labbel{cp}
Suppose that  $a $ is a $1$-compact element of  $\mathbf P$,
or just that $a$ is very-weakly-$1$-compact. 

Then $a$ is a non-generator if and only if 
 $a \leq p$, for every $p$ maximal$_{_{<1}}$  in $\mathbf P$.
 \end{corollary}

 \begin{proof}
Immediate from Proposition \ref{cp}(1) $\Leftrightarrow $  (4).  
 \end{proof}

If the meet of all the maximal$_{_{<1}}$  elements
of $\mathbf P$ exists, it shall be denoted
by $\Phi$ and shall be called the \emph{Frattini element}
of $\mathbf P$.
If $\mathbf P$ has no maximal$_{_{<1}}$  element,
we set $\Phi= 1$.

An element $b$ of $\mathbf P$ is \emph{locally $1$-compact}
(\emph{very-weakly-locally-$1$-compact})
if $b$ can be expressed as the join of $1$-compact 
(very-weakly-$1$-compact) elements.  
A poset $\mathbf P$ is (very-weakly)-$1$-algebraic if every element of $P$
is    (very-weakly)-locally $1$-compact.

\begin{corollary} \labbel{corcp}
Suppose that $\mathbf P$ is a poset with maximum $1$.
  \begin{enumerate}[(1)]
    \item
If $\Phi$ exists and is locally $1$-compact, or just very-weakly-$1$-compact, 
then $\Phi$ is the join of all the non-generators
of $\mathbf P$.
\item
In particular, if  $\mathbf P$ is a complete $1$-algebraic lattice,
or just a  complete very-weakly-$1$-algebraic lattice,
then  $\Phi$ is the join of all the non-generators
of $\mathbf P$.
\item
If $\mathbf P$ is a very-weakly-$1$-algebraic poset, then 
$\Phi$ exists if and only if the join of all 
the non-generators exists and, if this is the case,
they coincide. 
\item
If $\mathbf P$ is a finite
join-semilattice and there exists at least one non-generator,
then 
$\Phi$ exists and is the join of all 
the non-generators.
   \end{enumerate} 
 \end{corollary} 

\begin{proof}
(1) $\Phi$ is greater than  all the non-generators, by
Theorem \ref{cpp}(1) $\Rightarrow $  (i).
By assumption, 
 $\Phi$ can be expressed as the join of a family of 
very-weakly-$1$-compact elements,
which are thus non-generators, because of 
Corollary \ref{cp}.
Hence every $r \in P$
which contains all the non-generators
must contain  $\Phi$.

(2) is a special case of (1).

(3) If $\Phi$ exists, then the join of all 
the non-generators exists by (1), and
they coincide, again by (1). 
 
Suppose that the join  of all 
the non-generators exists, call it $\Gamma$.
Every maximal$_{_{<1}}$  element $p$ is greater than all the non-generators,
by Theorem \ref{cpp}(1) $\Rightarrow $  (i), hence $p \geq \Gamma $,
since $\Gamma$  is   the join  of all 
the non-generators.

Suppose that $r \leq p$,  for every maximal$_{_{<1}}$  element $p$.
Since $\mathbf P$ is  very-weakly-$1$-algebraic,
then $r= \bigvee A$, where $A$ is a family
of very-weakly-$1$-compact elements.
Since $a \leq r \leq p$, for every $a \in A$ and every maximal$_{_{<1}}$ 
element $p$, then every $a$ in $A$  is a non-generator, by Corollary \ref{cp}.
Hence $ a \leq \Gamma$, for every $a \in A$. In particular, 
$ r \leq \Gamma$, since  $r= \bigvee A$,

We have showed that 
$ \Gamma \leq p$ for every maximal$_{_{<1}}$ element $p$.
Moreover, if $r \leq p$,  for every maximal$_{_{<1}}$  element $p$,
then $ r \leq \Gamma$.
This shows that $\Gamma$ 
is the meet of all the maximal$_{_{<1}}$ elements,
thus $\Phi$ exists and  $\Phi = \Gamma$. 

(4) is immediate from (3).
Notice that $\Phi$ might not exist in a finite join-semilattice $\mathbf S$:
just let  $\mathbf S$
consist of just a maximum $1$ and at least two 
maximal$_{_{<1}}$ elements   and no more element.
 \end{proof}

\section{Further remarks} \labbel{fur}

As we mentioned in the introduction, 
we now see that dealing with 
closure posets or even closure spaces
does not improve the results in the preceding section.
Actually, we only assume that some poset $\mathbf R$  has 
a specified subset $P$ of  ``closed'' elements
and we show that it is no loss of generality to work
already in $\mathbf P$. 

Suppose that $\mathbf R$  is a poset with maximum
$1$, $P \subseteq R$ and   $1 \in P$.
Elements of $P$ shall be called \emph{closed}.
However, we make no specific assumption on 
the induced poset $\mathbf P$,
for example, we do not assume that
$\mathbf P$ is a semilattice.   
We say that some subset $Y \subseteq R$
\emph{$P$-generates $1$} if there is no $b \in P$ 
such that  $b < 1$ and  $y \leq b$, for every $y \in Y$. 

We say that $a \in R$ is a $P$-non-generator
in case that,  for every $e \in R$,
if $\{ a, e\}$  $P$-generates $1$, then 
$\{  e\}$  $P$-generates $1$. 

\begin{lemma} \labbel{bast}
Suppose that $\mathbf R$  is a poset with maximum
$1$ and $P \subseteq R$.
Let $\mathbf P$ be the induced poset on $P$.
 
If $a \in R$ and there exists a smallest element
$ c \in P$ such that $a \leq c$,
then $a$ is  a $P$-non-generator in $\mathbf R$ 
if and only if $c$ is a non-generator 
in $\mathbf P$. 
\end{lemma} 

\begin{proof}
Suppose that
$a$ is  a $P$-non-generator
and $c \vee d = 1$ in $\mathbf P$.
Since $c$  is the smallest closed containing $a$,
then there is no closed $<1$ containing both $a$ and $d$.   
Thus $\{ a, d\}$  $P$-generates $1$,
hence $\{ d\}$  $P$-generates $1$,
 since  $a$ is  a $P$-non-generator.
But $d$ is closed, hence $d=1$.  

Conversely, assume 
that $c$ is a non-generator 
in $\mathbf P$ and
$\{ a, b\}$  $P$-generates $1$.
Suppose by contradiction that 
$\{ b\}$ does not  $P$-generate $1$,
hence $b \leq d <1$, for some $d \in P$.  
Since $\{ a, b\}$  $P$-generates $1$
and $a \leq c$,  $b \leq d $,
then 
$\{ c, d\}$  $P$-generates $1$,
but this amounts to say that 
$c \vee d=1$ in $\mathbf P$,
since $c$ and $d$ are closed.    
Since $c$ is a non-generator 
in $\mathbf P$, then $d=1$, 
a contradiction.
 \end{proof}

\begin{remark} \labbel{non}   
The present abstract study is not intended to capture all
the significant aspects of Frattini's and Jacobson's theory.
At the heart of their radical theory there seems to be a deep and unexpected
connection between subalgebra- or module-like notions on the
``maximality side'' and congruence-like notions at the radical level.
 Whether the group and ring theory of Frattini and Jacobson
can be fully extended to a general universal-algebraic setting,
or even to a categorical one, seems to be still an
open problem. In this connection, see however \cite{KV,K}.
\end{remark}

\end{document}